# COUNTING NILPOTENT GALOIS EXTENSIONS

JÜRGEN KLÜNERS AND GUNTER MALLE

ABSTRACT. We obtain strong information on the asymptotic behaviour of the counting function for nilpotent Galois extensions with bounded discriminant of arbitrary number fields. This extends previous investigations for the case of abelian groups. In particular, the result confirms a conjecture by the second author on this function for arbitrary groups in the nilpotent case. We further prove compatibility of the conjecture with taking wreath products with the cyclic group of order 2 and give examples in degree up to 8.

## 1. INTRODUCTION

In this paper we are interested in the density of number fields with a given Galois group. Let $G \leq \mathfrak{S}_n$ be a finite transitive permutation group on $n$ points and let $k$ be a number field. We say that a finite extension $K/k$ has Galois group $G$ if the normal closure $\hat{K}$ of $K/k$ has Galois group isomorphic to $G$ and $K$ is the fixed field in $\hat{K}$ under a point stabilizer of $G$. By abuse of notation we will write $\text{Gal}(K/k) = G$ in this situation. We let

$$Z(k, G; x) := |\{K/k \mid \text{Gal}(K/k) = G, \ \mathcal{N}_{k/\mathbb{Q}}(d_{K/k}) \leq x\}|$$

be the number of field extensions of $k$ (inside a fixed algebraic closure $\bar{\mathbb{Q}}$ of $\mathbb{Q}$) of degree $n$ with Galois group permutation isomorphic to $G$ (as explained above) and norm of the discriminant $d_{K/k}$ bounded above by $x$. It is well known that the number of extensions of $k$ with bounded norm of the discriminant is finite, hence $Z(k, G; x)$ is finite for all $G$, $k$ and $x \geq 1$. Note that $Z(k, G; x)$ depends on the actual permutation representation $\phi: G \hookrightarrow \mathfrak{S}_n$; if we want to stress this dependence, then we write $Z(k, \phi(G); x)$.

For any finite subset $S \subset \mathbb{P}(k)$ of the set of prime divisors of $k$ we also consider extensions *unramified* in $S$, with corresponding counting functions

$$Z(k, S, G; x) := |\{K/k \mid \text{Gal}(K/k) = G, \ K/k \text{ unram. in } S, \ \mathcal{N}_{k/\mathbb{Q}}(d_{K/k}) \leq x\}|.$$

In [11] the second author put forward a conjecture on the asymptotic behaviour of $Z(k, S, G; x)$ for arbitrary $G$ as $x$ tends to infinity. In order to state it, we need to introduce a group-theoretic invariant of permutation groups: For a permutation $g$ (on a finite set $M$) we define its *index* as

$$\text{ind}(g) := |M| - \text{ the number of orbits of } g \text{ on } M.$$

Then, for a permutation group $G \neq 1$ we let

$$a(G) := \left(\min\{\text{ind}(g) \mid 1 \neq g \in G\}\right)^{-1},$$

and we define $a(1) := 0$ for the trivial group.

**Conjecture 1.1.** [11] *Let $G$ be a transitive permutation group, $k$ a number field, and $S \subset \mathbb{P}(k)$ a finite set of primes of $k$. Then there exists a constant $c_1(k, S, G) > 0$ and for all $\epsilon > 0$ there exists a constant $c_2(k, G, \epsilon) > 0$ such that*

$$c_1(k, S, G) \, x^{a(G)} \leq Z(k, S, G; x) < c_2(k, G, \epsilon) \, x^{a(G)+\epsilon}$$

*for all large enough $x$.*

Our main result is Theorem 6.3, proving that Conjecture 1.1 holds for nilpotent groups in their regular permutation representation, that is, for Galois extensions of arbitrary number fields with nilpotent Galois group. One important ingredient (for the lower bound) is the celebrated Theorem of Šafarevič asserting that there exists at least one extension for any nilpotent group. Surprisingly enough, this is sufficient to deduce the correct lower bound.

We also prove the upper bound for arbitrary permutation representations of $\ell$-groups (Theorem 7.3). Moreover, we show that Conjecture 1.1 is true for wreath products of the cyclic group $C_2$ with an arbitrary group $H$, provided it holds for $H$ (Corollary 8.4), using a general argument from [11], and the lower bound holds for wreath products of cyclic groups of prime order with an arbitrary group $H$ if there exists at least one $H$-extension (Theorem 8.2).

An even stronger form of the conjecture would also control the constant $c_2$ in the upper bound as an explicit function of invariants of $k$. Such a form would be much better adapted to inductive proofs, see [11], Prop. 5.2.

We will sometimes consider the weaker form of Conjecture 1.1 obtained by choosing $S = \emptyset$. Note that as far as the upper bound is concerned this is no restriction.

Previously, the truth of Conjecture 1.1 was only known to hold for abelian groups by work of Wright [17] (see also [4]), and for the non-abelian groups $\mathfrak{S}_3$ by Davenport and Heilbronn, and $D_4$ by Cohen, Diaz y Diaz and Olivier [3] (in all cases in a much more precise form). The lower bound for $\mathfrak{A}_4$ and for $\mathfrak{S}_4$ over $k = \mathbb{Q}$ was shown by Baily [1]. We prove it in the case of arbitrary $k$ in Proposition 10.2.

**Acknowledgement:** We would like to thank Florin Nicolae for the carefully reading of a preliminary version of this work.

## 2. Notation and auxiliary results

In this section we prove some lemmata which are useful to get the desired bounds. All number fields will be considered lying inside a fixed algebraic closure of $\mathbb{Q}$. Throughout, for a number field $k$ we write $d_k$ for the discriminant $d_{k/\mathbb{Q}}$. Also, the absolute norm of $k/\mathbb{Q}$ will be denoted by $\mathcal{N}$.

For brevity we will adopt the notation 'for all $x \gg 0$' to mean 'for all sufficiently large $x > 0$'.

Let $k_1$ be a number field and $S \subset \mathbb{P}(k_1)$ a set of primes. Then for $b \in \mathbb{Z}$ we write $b^S$ for the part of $b$ which is coprime to the norms of elements in $S$. For a number field $k_2$ we let $Y^S(k_2, G; x)$ be the number of $G$-extensions of $k_2$ with $\mathcal{N}(d_{K/k_2})^S \leq x$. We get the following:

**Lemma 2.1.** *Let $G$ be a finite group, $S$ a finite set of primes and $m \geq 1$. Then there exists a constant $c(m, S, G) \in \mathbb{R}$ only depending on $m$, $S$ and $G$ such that, whenever $k$ is a number field of degree $[k : \mathbb{Q}] = m$ with $c(k, S) x^a \leq Y^S(k, G; x)$ for some $c(k, S) \in \mathbb{R}$, $a \geq 0$ and all $x \gg 0$, then we have*

$$c(k, S) c(m, S, G)^{-a} x^a \leq Z(k, G; x) \leq Y^S(k, G; x)$$

*for all $x \gg 0$.*

*Proof.* By definition we always have $Z(k, G; x) \leq Y^S(k, G; x)$. Let $p \in S$ be a prime. There is a well-known upper bound for the maximal possible $p$-power $p^e$ dividing $\mathcal{N}(d_{K/k})$, where $K$ is a $G$-extension of $k$ (see [16], Remark 1, p. 58). This bound only depends on $[K : \mathbb{Q}]$ and $p$, hence on $[k : \mathbb{Q}] = m$, $G$ and $p$. The product over all these bounds for the primes $p \in S$ gives an integer $c(m, S, G)$ such that $c(m, S, G) \mathcal{N}(d_{K/k})^S \geq \mathcal{N}(d_{K/k})$ for all $G$-extensions $K/k$. This implies

$$Z(k, G; x) \geq Y^S(k, G; x/c(m, S, G)) \geq c(k, S)/c(m, S, G)^a \, x^a$$

for all $x \gg 0$ as claimed. □

The above lemma enables us to neglect the contribution of a finite number of primes when we want to determine a lower bound for the number of $G$-extensions of $k$. Usually we will choose $S$ such that it contains the primes smaller than $[K:k]$. This already avoids having to consider wild ramification.

Sometimes the set $S$ which we want to use changes. Let $\mathfrak{a}$ be an ideal in the ring of integers $\mathcal{O}_k$ of a number field $k$. We write $t_k(\mathfrak{a})$ for the number of ideal divisors of $\mathfrak{a}$ in $\mathcal{O}_k$. We need the following estimate.

**Lemma 2.2.** *For all $\epsilon > 0, m \in \mathbb{N}$, there exists a constant $c = c(m, \epsilon)$ such that, whenever $k$ is a number field of degree $[k:\mathbb{Q}] = m$, then $t_k(\mathfrak{a}) \leq c\mathcal{N}(\mathfrak{a})^\epsilon$ for all $\mathfrak{a} \triangleleft \mathcal{O}_k$.*

*Proof.* We imitate the proof for the case $k = \mathbb{Q}$ in [8], Thm. 315. Clearly $t_k(\mathfrak{a}) = \prod_\mathfrak{p}(e_\mathfrak{p} + 1)$ when $\mathfrak{a} = \prod_\mathfrak{p} \mathfrak{p}^{e_\mathfrak{p}}$ is the prime ideal factorization of $\mathfrak{a} \triangleleft \mathcal{O}_k$. Hence

$$\frac{t_k(\mathfrak{a})}{\mathcal{N}(\mathfrak{a})^\epsilon} = \prod_\mathfrak{p} \frac{e_\mathfrak{p} + 1}{\mathcal{N}(\mathfrak{p})^{e_\mathfrak{p}\epsilon}}.$$

We may estimate the factors by

$$\frac{e_\mathfrak{p} + 1}{\mathcal{N}(\mathfrak{p})^{e_\mathfrak{p}\epsilon}} \leq 1 + \frac{e_\mathfrak{p}}{\mathcal{N}(\mathfrak{p})^{e_\mathfrak{p}\epsilon}} \leq 1 + \frac{1}{\epsilon \log 2} \leq \exp(\frac{1}{\epsilon \log 2})$$

since $c \log 2 \leq e^{c \log 2} = 2^c \leq s^c$ for all $c > 0$, $s \geq 2$. Moreover, for prime ideals $\mathfrak{p}$ with $\mathcal{N}(\mathfrak{p})^\epsilon \geq 2$ we have

$$\frac{a+1}{\mathcal{N}(\mathfrak{p})^{a\epsilon}} \leq \frac{a+1}{2^a} \leq 1.$$

There exist no more than $2^{1/\epsilon}$ primes below $2^{1/\epsilon}$, hence no more than $[k:\mathbb{Q}]2^{1/\epsilon}$ prime ideals $\mathfrak{p}$ in $\mathcal{O}_k$ with $\mathcal{N}(\mathfrak{p}) \leq 2^{1/\epsilon}$. So

$$\frac{t_k(\mathfrak{a})}{\mathcal{N}(\mathfrak{a})^\epsilon} \leq \prod_{\mathcal{N}(\mathfrak{p})^\epsilon < 2} \exp(\frac{1}{\epsilon \log 2}) \leq \exp\left(\frac{[k:\mathbb{Q}]2^{1/\epsilon}}{\epsilon \log 2}\right),$$

a constant depending only on $m$ and $\epsilon$. □

Let $k$ be a number field and $I_k$ the group of fractional ideals of $\mathcal{O}_k$. Furthermore let $S$ be a finite set of prime ideals in $\mathcal{O}_k$. For any $\alpha \in \mathcal{O}_k$ we can decompose the principal ideal $(\alpha)$ as

$$(\alpha) = \prod_{\mathfrak{p} \in S} \mathfrak{p}^{e_\mathfrak{p}} \prod_{\mathfrak{p} \notin S} \mathfrak{p}^{e_\mathfrak{p}},$$

where both products are finite. For any $\ell \in \mathbb{P}$ the mapping defined by

$$\tilde{\phi} : \mathcal{O}_k \longrightarrow I_k, \qquad \alpha \mapsto \tilde{\phi}(\alpha) := \prod_{\mathfrak{p} \notin S} \mathfrak{p}^{e_\mathfrak{p}},$$

induces a map $\phi : \mathcal{O}_k^\times/(\mathcal{O}_k^\times)^\ell \to I_k/I_k^\ell$ modulo $\ell$th powers. We get the following lemma:

**Lemma 2.3.** *Let $k$ be a number field, $S \subset \mathbb{P}(k)$ finite and $\ell \in \mathbb{P}$. Then, all fibers of the map $\phi$ described above are finite of size at most $\ell^{|S|+r+u+1}$, where $u$ denotes the unit rank of $\mathcal{O}_k$ and $r$ the $\ell$-rank of the class group of $k$.*

*Proof.* First consider the natural map $\psi : \mathcal{O}_k^\times \to I_k/I_k^\ell$ sending an element to its principal ideal modulo $\ell$th powers. Let $B$ be a system of representatives of the class group in $I_k$. We claim that the kernel of $\psi$ is generated by $(\mathcal{O}_k^\times)^\ell$, the unit group and the elements of order $\ell$ in $B$. Indeed, if $\alpha \in \ker(\psi)$ then $(\alpha) := \alpha\mathcal{O}_k = \mathfrak{a}^\ell$ for some ideal $\mathfrak{a} \in I_k$. Let $\mathfrak{b} \in B$ be the representative of $\mathfrak{a}$, such that $\mathfrak{a}\mathfrak{b}^{-1} = (\beta)$ is principal. Then we have $(\alpha) = \mathfrak{a}^\ell = (\beta)^\ell \mathfrak{b}^\ell = (\beta^\ell)\mathfrak{b}^\ell$. Hence $\mathfrak{b}^\ell = (\gamma)$ is principal, so $\mathfrak{b}$ has order dividing $\ell$ in the class group. We obtain that $(\alpha) = (\beta^\ell \gamma)$, which proves that $\alpha$ and $\beta^\ell \gamma$ differ by a unit.

Thus, the kernel of $\phi$ is generated by units modulo $\ell$th powers, the elements whose principal ideal is supported by $S$ modulo $\ell$th powers, and the elements of order $\ell$ in $B$. It is then clear it has the desired size. $\square$

Throughout we will use the following well-known convergence criterion for Dirichlet series (see for example [15], Lemma 58.B):

**Lemma 2.4.** *Let $f(s) := \sum_{n \geq 1} \frac{a_n}{n^s}$ be a Dirichlet series with non-negative coefficients $a_n$. If $\sum_{n=1}^{x} a_n \leq c\, x^r$ for all $x$, then $f(s)$ converges for all $s > r$.*

Recall the counting function $Z(k, G; x)$ from the introduction for $G$-extensions of a number field $k$ with norm of the discriminant bounded by $x$. Let $\ell$ be a prime. Wright [17] has shown the upper bound $Z(k, C_\ell; x) \leq c(k, \ell)\, \sqrt[\ell-1]{x} \log^{\ell-2} x$ for $x \gg 0$, for the growth of the number of $C_\ell$-extensions of any number field $k$ (see also [4] for a more precise version). Here, we will prove an explicit form for the constant $c(k, \ell)$, only depending on the degree $[k : \mathbb{Q}]$ and the $\ell$-rank of the class group of $k$, which does not seem to follow from [17] in any obvious way. Let $\zeta$ be a primitive $\ell$th root of unity and write $\tilde{k} := k(\zeta)$.

**Theorem 2.5.** *Let $\ell$ be a prime and $m \geq 1$. Then for any $\epsilon > 0$ there exists a constant $c(m, \ell, \epsilon)$ such that for all number fields $k$ of degree $[k : \mathbb{Q}] = m$ and all finite subsets $S \subset \mathbb{P}(k)$ and $\mathfrak{q} \in \mathbb{P}(k)$ we have*

$$Y^S(k, C_\ell; x) \leq c(m, \ell, \epsilon)\, |\mathrm{Cl}(\tilde{k})_\ell|\, D^{r\epsilon}\, x^{a(C_\ell)+\epsilon} \qquad \text{for all } x > 0,$$

*where $D = \prod_{\mathfrak{p} \in S, \mathfrak{p} \neq \mathfrak{q}, \mathfrak{p} \nmid (\ell)} \mathcal{N}(\mathfrak{p})$, $\mathrm{Cl}(\tilde{k})_\ell$ denotes the largest elementary abelian $\ell$-factor group of the class group $\mathrm{Cl}(\tilde{k})$ and $r \geq 1$ is a constant only depending on $\ell$ and $m$.*

*In particular the case $S = \emptyset$, $D = 1$ gives*

$$Z(k, C_\ell; x) \leq c(m, \ell, \epsilon)\, |\mathrm{Cl}(\tilde{k})_\ell|\, x^{a(C_\ell)+\epsilon} \qquad \text{for all } x > 0.$$

*Proof.* Let $k/\mathbb{Q}$ be of degree $m$. By Kummer theory, $C_\ell$-extensions of $k$ are parametrized by elements $\alpha \in \tilde{k}^\times/(\tilde{k}^\times)^\ell$ (satisfying some additional condition, see Section 3). Denote the corresponding extension by $k_\alpha$.

Let $D' := \prod_{\mathfrak{p} \in S, \mathfrak{p} \mid (\ell), \mathfrak{p} \neq \mathfrak{q}} \mathcal{N}(\mathfrak{p})$. Then clearly $D'$ is bounded above by a constant depending on $\ell, m$ only. The function $Y^S(k, C_\ell; x)$ counts cyclic extensions $K/k$ with $\mathcal{N}(d_{K/k})^S = d_1$, where $d_1 \leq x$. That is, we count cyclic extensions $K/k$ with $\mathcal{N}(d_{K/k})$ dividing $d_2 = d_1(DD'\mathcal{N}(\mathfrak{q}))^r$, where $d_1 \leq x$ and the exponent $r$ only depends on $\ell$ and $m$ (see [16], Remark 1, p.58, for example).

Let's fix $d_2 \in \mathbb{N}$ and count cyclic extensions $k_\alpha$ of $k$ with norm of the discriminant $d = \mathcal{N}(d_{K/k})$ dividing $d_2$. Note that any such $d$ is $\ell-1$-powerful, that is, any prime divisor of $d$ occurs (at least) to the power $\ell - 1$. All prime divisors of the principal ideal $\alpha \mathcal{O}_{\tilde{k}}$ which do not ramify in $k_\alpha/k$ occur to a power divisible by $\ell$. So up to an $\ell$th power $\alpha \mathcal{O}_{\tilde{k}}$ divides $d \mathcal{O}_{\tilde{k}}$, hence $d_2 \mathcal{O}_{\tilde{k}}$, in $I_{\tilde{k}}$. Since $[\tilde{k} : \mathbb{Q}] \leq m(\ell-1)$, for any $\delta > 0$ there is a constant $c(m, \delta) > 0$ such that there exist at most $c(m, \delta) \mathcal{N}(d_1(DD')^r \mathcal{O}_{\tilde{k}})^\delta$ such divisors which are coprime to $\mathfrak{q}\, \mathcal{O}_{\tilde{k}}$ by Lemma 2.2. Altogether there are at most $(\ell-1)^r \mathcal{N}(d_1(DD')^r \mathcal{O}_{\tilde{k}})^\delta$ such divisors. Application of Lemma 2.3 with $S = \emptyset$ shows that for any principal ideal there exist at most $\ell^m |\mathrm{Cl}(\tilde{k})_\ell|$ generators $\alpha \in \tilde{k}$ modulo $\ell$th powers. But $\mathcal{N}_{\tilde{k}/\mathbb{Q}}((d_1(DD')^r \mathcal{O}_{\tilde{k}})^\delta \leq (d_1(DD')^r)^{m(\ell-1)\delta}$, so we find at most $\ell^m |\mathrm{Cl}(\tilde{k})_\ell| \tilde{c}(m, \ell, \delta)\, (d_1(DD')^r)^{m(\ell-1)\delta}$ possible $k_\alpha$ for a given $d_2 = d_1(DD'\mathcal{N}(\mathfrak{q}))^r$.

Let's define $\chi$ to be the characteristic function of the set of $\ell-1$-powerful integers. Hence

$$\begin{aligned} Y^S(k,C_\ell;x) &\leq \ell^m |\mathrm{Cl}(\tilde{k})_\ell| \tilde{c}(m,\ell,\delta) \sum_{d_1=1}^{x} \chi(d_1) \left(d_1 (DD')^r\right)^{m(\ell-1)\delta} \\ &\leq \ell^m |\mathrm{Cl}(\tilde{k})_\ell| \tilde{c}(m,\ell,\delta) (DD')^{rm(\ell-1)\delta} x^{\frac{1}{\ell-1}+\epsilon} \sum_{d=1}^{x} \frac{\chi(d)}{d^{\frac{1}{\ell-1}+\epsilon-m(\ell-1)\delta}}. \end{aligned}$$

But by a result of Erdös and Szekeres [6], the number of $\ell-1$-powerful integers below $x$ grows asymptotically like $c_1 x^{\frac{1}{\ell-1}}$ for some constant $c_1 > 0$. Hence we reach the desired conclusion from Lemma 2.4 by choosing $m(\ell-1)\delta < \epsilon$ (with constants only depending on $\epsilon$). □

We also need a lower bound for the number of cyclic extensions of prime degree unramified in a finite set $S$. This follows easily from Wright's result for the case $S = \emptyset$:

**Theorem 2.6.** *Let $\ell$ be a prime and $k$ a number field. Then for any finite subset $S \subset \mathbb{P}(k)$ there exists a constant $c(\ell, S) > 0$ such that*

$$c(\ell, S) x^{a(C_\ell)} \leq Z(k, S, C_\ell; x) \qquad \text{for all } x \gg 0.$$

*Proof.* By the main theorem of [17] there exists a constant $c > 0$ such that $c x^{a(C_\ell)} \leq Z(k, C_\ell; x)$, so the statement holds for $|S| = 0$. Now assume that $|S| > 0$ and let $\mathfrak{p} \in S$. By induction we may assume that there exists $c(S') > 0$ with $c(S') x^{a(C_\ell)} \leq Z(k, S', C_\ell; x)$ for $x \gg 0$, where $S' := S \setminus \{\mathfrak{p}\}$. Consider sets $\{L_1, \ldots, L_r\}$ of $C_\ell$-extensions of $k$ unramified in $S'$ such that every $C_\ell$-subextension of the composite $L_1 \cdots L_r/k$ is ramified in $\mathfrak{p}$. Since the completion $k_\mathfrak{p}$ has only finitely many extensions of degree $\ell$, $r$ is bounded from above.

Choose $L_1, \ldots, L_r$ as above with $r$ maximal and let $L$ be their composite. Then for any $C_\ell$-extension $K/k$ not contained in $L$ the composite $KL/k$ contains a $C_\ell$-extension $K_u/k$ unramified in $\mathfrak{p}$. In particular, if $K/k$ is unramified in $S'$, then $K_u/k$ is unramified in $S$. Moreover, if $K'$ is different from the $s := r\ell + 1$ subfields of $KL/k$ with group $C_\ell$, then the corresponding extension $K'_u$ unramified in $\mathfrak{p}$ is different from $K_u$. Thus for any $s$ distinct $C_\ell$-extensions of $k$ we obtain at least one $C_\ell$-extension unramified in $S$. Finally, by the discriminant composition formula the norm of $d_{K_u/k}/d_{K/k}$ is bounded above by a constant only depending on $k$ and $\ell$. Thus there exists a constant $c(S) > 0$ with $c(S) x^{a(C_\ell)} \leq Z(k, S, C_\ell; x)$ for all $x \gg 0$. □

We may even choose the decomposition types at finitely many places:

**Corollary 2.7.** *Theorem 2.6 remains true if we prescribe the decomposition types at the finitely many (unramified) places in $S$.*

*Proof.* The same argument as in the previous proof gives the desired result by induction on $|S|$, using that the composite of any two distinct $C_\ell$-extensions unramified at $\mathfrak{p}$ contains both a completely split and an inert $C_\ell$-extension. □

The following lemma, due to Šafarevič, will be used to investigate solutions to embedding problems:

**Lemma 2.8.** *Let $k$ be a field of characteristic different from $\ell$ containing an $\ell$th root of unity $\zeta$ and $\hat{k} := k(\sqrt[\ell]{\alpha})$ for some $\alpha \in k^\times \setminus (k^\times)^\ell$. Let $\sigma \in \mathrm{Aut}(k)$ such that $\hat{k}/k^{\langle\sigma\rangle}$ is Galois. Then $\hat{k}/k^{\langle\sigma\rangle}$ is abelian if and only if $\alpha^\sigma = \alpha^q \beta^\ell$ for some $\beta \in k$, where $q \in \mathbb{N}$ is defined by $\zeta^\sigma = \zeta^q$.*

*Proof.* Since $\hat{k}/k^{\langle\sigma\rangle}$ is Galois we have $k(\sqrt[\ell]{\alpha}) = k(\sqrt[\ell]{\alpha^\sigma})$. By [13], Lemma 7.15, this implies that $\alpha^\sigma = \alpha^r \beta^\ell$ for some $\beta \in k^\times$ and some $1 \leq r \leq \ell - 1$. Let $\rho$ denote a generator of $\mathrm{Gal}(\hat{k}/k)$, with $\rho : \sqrt[\ell]{\alpha} \mapsto \zeta\sqrt[\ell]{\alpha}$. Lift $\sigma$ to an automorphism $\hat{\sigma}$ of $\hat{k}/k^{\langle\sigma\rangle}$, so that
$$\hat{\sigma} : \sqrt[\ell]{\alpha} \mapsto \sqrt[\ell]{\alpha^r}\beta, \quad \zeta \mapsto \zeta^q.$$
Then $\rho\hat{\sigma} = \hat{\sigma}\rho$ if and only if $q \equiv r \pmod{\ell}$, as claimed. $\square$

Most of our results deal with the regular permutation representation of a finite group $G$. Note that the value of $a(G)$ can easily be given in that case:

**Remark 2.9.** *Let $G \leq \mathfrak{S}_n$ be a regular permutation group. Denote by $\ell$ the smallest prime dividing $|G|$. Then $a(G) = \frac{\ell}{(\ell-1)|G|} = \frac{\ell}{(\ell-1)n}$.*

## 3. Solutions to central embedding problems

In this section we prepare the proof of Conjecture 1.1 for normal extensions with nilpotent Galois group of an arbitrary number field $k$. For this, we study so-called embedding problems.

**Definition 3.1.** Let $K/k$ be a Galois extension of number fields with group $H$ and
$$1 \longrightarrow U \longrightarrow G \stackrel{\kappa}{\longrightarrow} H \longrightarrow 1$$
an exact sequence of finite groups. Denote by $\bar{k}$ an algebraic closure of $k$ and by $\varphi : \mathrm{Gal}(\bar{k}/k) \to H \cong \mathrm{Gal}(K/k)$ an epimorphism with kernel $\mathrm{Gal}(\bar{k}/K)$. This is called the *embedding problem given by $\kappa$ and $\varphi$*. An epimorphism $\tilde{\varphi} : \mathrm{Gal}(\bar{k}/k) \to G$ is called *(proper) solution* to the embedding problem, if $\kappa \circ \tilde{\varphi} = \varphi$. The fixed field $L$ of $\ker(\tilde{\varphi})$ is called *solution field*. Furthermore the embedding problem is called *central*, if $U$ is in the center of $G$.

Suppose that we are given an extension $K/k$ with Galois group $H$ and we want to find all fields $L$ containing $K$ such that $\mathrm{Gal}(L/k) \cong G$. Since the embedding problems depend on $\varphi$ and $\kappa$ these can in general not be obtained as solutions of a single embedding problem.

**Lemma 3.2.** *Let $K/k$ be an extension of number fields with Galois group $H$ and let $G$ be a finite central extension of $H$. Then all fields $L/K$ such that $\mathrm{Gal}(L/k) \cong G$ can be found as solutions of finitely many embedding problems, whose number only depends on $G$ and $H$.*

*Proof.* Let $\varphi$ be defined as in Definition 3.1. Epimorphisms $\kappa : G \to H$ with the same kernel only differ by automorphisms of $H$. Furthermore $G$ has only finitely many normal subgroups with quotient isomorphic to $H$. $\square$

Let first $H$ be an arbitrary finite group and
$$(*) \qquad\qquad 1 \longrightarrow C_\ell \longrightarrow G \longrightarrow H \longrightarrow 1$$
a central embedding problem with kernel $C_\ell$. Later we specialize to the case that $H$ is an $\ell$-group or nilpotent. Throughout the next four sections, we will be concerned with regular permutation representations only, that is, we count *Galois* extensions with given Galois group. In Corollary 7.3 we obtain from this a result valid for arbitrary permutation representations of $\ell$-groups.

We want to study the asymptotics of the number of $G$-extensions of $k$, assuming some knowledge on the asymptotics of $H$-extensions. If the ground field contains the $\ell$th roots of unity, (*) becomes a Brauer embedding problem, for which the question of solvability of (*) and of a parametrization of all solutions is much easier. So we will also have to consider that case. Throughout, we fix a primitive $\ell$th root of unity $\zeta$ and we write $\tilde{k} := k(\zeta)$ for any number field $k$. We set $z := [\tilde{k} : k]$.

We now characterize solution fields to the embedding problem for $K/k$. In order not to obscure the presentation we first prove the version without restriction on the ramification.

**Proposition 3.3.** *Let $K/k$ be a Galois extension with group $H$ and $L$ a solution field to the embedding problem (\*) for $K/k$. Let $\sigma$ be a generator of $\mathrm{Gal}(\tilde{L}/L)$ with $\zeta^\sigma = \zeta^q$. Then there exists $\alpha \in \tilde{K}^\times$ such that $\tilde{L} = \tilde{K}(\sqrt[\ell]{\alpha})$ and $\alpha^\sigma = \alpha^q \beta^\ell$ for some $\beta \in \tilde{K}$.*

*Proof.* Let $U \leq H$ be the Galois group of $K$ over $K \cap \tilde{k}$, and $s := [\tilde{K} : K]$, so $(H : U) = z/s$. Then $\tilde{L}/k$ is Galois with group a subdirect product $G \times_{z/s} C_z$, and $\tilde{L}$ is a solution field to the induced embedding problem
$$1 \longrightarrow C_\ell \longrightarrow C_\ell.U \longrightarrow U \longrightarrow 1$$
for $\tilde{K}/\tilde{k}$. By Kummer theory there exists $\alpha \in \tilde{k}^\times$ with $\tilde{L} = \tilde{K}(\sqrt[\ell]{\alpha})$. Since $\tilde{L}/K$ is abelian by construction, the assertion follows from Lemma 2.8 with $k = \tilde{K}$, $\hat{k} = \tilde{L}$. □

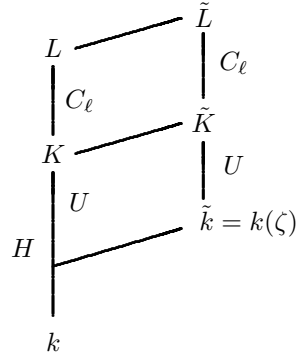

Now assume fixed a solution $L/k$ of the embedding problem (\*) for $K/k$, with $\sigma, \alpha, q$ as in Proposition 3.3. Let $L'$ be another solution field to (\*). As in the proof of Proposition 3.3 we see that $\tilde{L}'$ is a solution field to the induced embedding problem for $\tilde{K}/\tilde{k}$. The latter is a Brauer embedding problem, since the ground field contains the $\ell$th roots of unity. By [12], Thm. IV.7.2, any such solution field is of the form $\tilde{K}(\sqrt[\ell]{b\alpha})$ for some $b \in \tilde{k}^\times$. Since $\tilde{L}' = L'\tilde{K}$ and $L'/K$ and $\tilde{K}/K$ are abelian, the extension $\tilde{L}'/K$ is abelian. By Lemma 2.8 the element $(b\alpha)^\sigma/(b\alpha)^q$ is an $\ell$th power in $\tilde{K}$, so the same is true for $b^\sigma/b^q$ by the property of $\alpha$.

We claim that even $b^\sigma/b^q \in \tilde{k}^\ell$. Indeed, the extension $\tilde{k}(\sqrt[\ell]{b})/k$ is Galois. But then its translation $\tilde{K}(\sqrt[\ell]{b})/K$ with $K/k$ is Galois with the same group. Since $b^\sigma/b^q \in \tilde{K}^\ell$ this latter extension is abelian by Lemma 2.8, so $\tilde{k}(\sqrt[\ell]{b})/k$ is abelian. Again by Lemma 2.8 this proves that $b^\sigma/b^q \in \tilde{k}^\ell$.

For the converse, assume that $K$ is linearly disjoint from $\tilde{k}$ over $k$. Let $b \in \tilde{k}^\times$ with $b^\sigma/b^q \in \tilde{k}^\ell$. Then $\tilde{L}_b := \tilde{K}(\sqrt[\ell]{b\alpha})$ is a solution field to the induced embedding problem over $\tilde{k}$. By [13], Lemma 7.15, and the fact that $(b\alpha)^\sigma/(b\alpha)^q \in \tilde{K}^\ell$ the extension $\tilde{L}_b/k$ is Galois, with Galois group an extension of $\mathrm{Gal}(\tilde{K}/k) = H \times C_z$ by $C_\ell$. By Lemma 2.8 the extension $\tilde{L}_b/K$ is abelian, hence has group $C_\ell \times C_z$. The fixed field $L_b$ under the unique normal subgroup $C_z$ of $\mathrm{Gal}(\tilde{L}_b/K)$ thus has group $G$ over $k$, hence is a solution to the original embedding problem.

We have thus proved the following classification of all solution fields:

**Proposition 3.4.** *Let $K/k$ be a Galois extension with group $H$, linearly disjoint from $\tilde{k}/k$, $L$ a solution field to the embedding problem (\*) for $K/k$, and $\alpha, \sigma, q$ as in Proposition 3.3.*

*Then the solutions to the embedding problem (\*) for $K/k$ are precisely the fixed fields inside $\tilde{L}_b := \tilde{K}(\sqrt[\ell]{b\alpha})$ under the subgroup of order $z = [\tilde{k} : k]$ of $\mathrm{Gal}(\tilde{L}_b/K) \cong C_\ell \times C_z$, for elements $b \in \tilde{k}^\times$ satisfying $b^\sigma/b^q \in (\tilde{k}^\times)^\ell$.*

In the case that $K/k$ is not linearly disjoint from $\tilde{k}/k$, the induced embedding problem for $\tilde{K}/\tilde{k}$ involves a strictly smaller group. It is not clear which solutions to this induced problem lead to solutions of the original embedding problem for $K/k$.

The preceding result yields a very strong connection between cyclic extensions of $k$ of degree $\ell$ and solutions to the original embedding problem for $K/k$. Indeed, specializing to the case $H = 1$, $K = k$, we see that both are indexed by the same set of parameters $b$. Moreover, this connection has finite bounded fibers:

**Proposition 3.5.** *Let $k$ be a number field of degree $m = [k : \mathbb{Q}]$, $K/k$ a Galois extension with group $H$ linearly disjoint from $\tilde{k}/k$, and assume that the embedding problem (\*) for $K/k$ is solvable. Then we have:*

(a) *The correspondence $k_b \mapsto L_b$ between cyclic Galois extensions of $k$ with group $C_\ell$ and solutions to the embedding problem (\*) for $K/k$ has finite fibers of size bounded only in terms of $H$, $\ell$ and $m$.*

(b) *Moreover, the discriminants satisfy*

$$c_1(K)\mathcal{N}(d_{L_b/k}) \le \mathcal{N}(d_{k_b/k})^{|H|} \le c_2(K)\mathcal{N}(d_{L_b/k})$$

*for constants $c_1(K), c_2(K) > 0$ only depending on $K$ (and $\ell$).*

*Proof.* By Proposition 3.4 any cyclic extension of degree $\ell$ of $k$ is of the form $k_b$ for a unique $b \in \tilde{k}^\times/(\tilde{k}^\times)^\ell$ with $b^\sigma/b^q \in (\tilde{k}^\times)^\ell$. This determines a unique field $\tilde{L}_b = \tilde{K}(\sqrt[\ell]{b\alpha})$, which in turn corresponds to a unique solution field $L_b$ to (\*).

Conversely, let $\alpha, \sigma, q$ be as in Proposition 3.3. Assume given $b_i$, $1 \le i \le t$, such that $L_{b_i}$ are mutually isomorphic. Then so are the fields $\tilde{L}_i := \tilde{L}_{b_i} = \tilde{K}(\sqrt[\ell]{b_i\alpha})$, $1 \le i \le t$. In particular $\sqrt[\ell]{b_1/b_j} \in \tilde{L}_1$ for all $j$, hence $\tilde{k}(\sqrt[\ell]{b_1/b_j})$ is a subfield of $\tilde{L}_1/\tilde{k}$. Since the latter has only finitely many subfields, for large enough $t$ (only depending on $G$) some of the $\tilde{k}(\sqrt[\ell]{b_1/b_j})$ have to agree, which in turn forces some of the $k_{b_i}$ to coincide. This completes the proof of assertion (a).

By the composition formula the discriminants of the various fields are related by

$$d_{\tilde{L}_b/k} = d_{\tilde{k}/k}^{|G|}\mathcal{N}_{\tilde{k}/k}(d_{\tilde{L}_b/\tilde{k}}) = (d_{L_b/k})^z \mathcal{N}_{L_b/k}(d_{\tilde{L}_b/L_b})$$

with $z := [\tilde{k} : k]$. Clearly, the norms of both $d_{\tilde{k}/k}^{|G|}$ and $\mathcal{N}_{L_b/k}(d_{\tilde{L}_b/L_b})$ are non-zero integers bounded independently from $b$ only in terms of $\ell$ and $[k : \mathbb{Q}] = m$. So $\mathcal{N}(d_{L_b/k})^z$ only differs by bounded factors only depending on $\ell, m$ from $\mathcal{N}(d_{\tilde{L}_b/\tilde{k}})$. In particular, choosing $H = 1$, $K = k$, we conclude that $\mathcal{N}(d_{k_b/k})^z$ and $\mathcal{N}(d_{\tilde{k}_b/\tilde{k}})$ only differ by globally bounded factors. Thus it suffices to prove (b) in the case when $k = \tilde{k}$ contains the $\ell$th roots of unity.

Application of the discriminant composition formula to the extension $M/k := L_b k_b/k$ yields

$$d_{L_b/k}^\ell \mathcal{N}_{L_b/k}(d_{M/L_b}) = d_{k_b/k}^{|G|}\mathcal{N}_{k_b/k}(d_{M/k_b}) = d_{k_b/k}^{|G|}\mathcal{N}_{k_b/k}(d_{K_b/k_b})^\ell \mathcal{N}_{K_b/k}(d_{M/K_b})$$

with $K_b := Kk_b$. Now both extensions $M/L_b$ and $M/K_b$ are generated by $\sqrt[\ell]{\alpha}$, hence their discriminants only contribute bounded terms. Furthermore, the norm of $\mathcal{N}_{k_b/k}(d_{K_b/k_b})$ is bounded by the norm of the discriminant of $K/k$. We conclude that $\mathcal{N}(d_{L_b/k})^\ell$ and $\mathcal{N}(d_{k_b/k})^{|G|}$ only differ by globally bounded factors, as claimed. □

In order to treat the case of controlled ramification we need the following result of Hecke (see [2], Thm. 10.2.9, for example):

**Theorem 3.6.** (Hecke) *Let $k$ be a number field containing the $\ell$th roots of unity and $N = k(\sqrt[\ell]{\alpha})$. Let $\mathfrak{p}$ be a prime of $\mathcal{O}_k$. If $\mathfrak{p}$ is ramified in $N/k$ then either $\nu_\mathfrak{p}(\alpha) \not\equiv 0 \pmod{\ell}$ or $\mathfrak{p}|\ell$. If $\mathfrak{p}|\ell$ and $\alpha$ is prime to $\ell$, then $\mathfrak{p}$ is unramified in $N/k$ if the congruence $x^\ell \equiv \alpha$ has a solution in $k$ modulo $\mathfrak{p}^a$ for some $a \ge (\ell-1)\nu_\mathfrak{p}(1-\zeta)$.*

This allows to choose solutions to (\*) with ramification outside a given set as follows:

**Lemma 3.7.** *Let $k$ be a number field containing the $\ell$th roots of unity and $S \subset \mathbb{P}(k)$ be finite. Furthermore let $1 \to C_\ell \to G \to H \to 1$ be a central group extension, $K/k$ a normal extension with Galois group $H$ ramified in $S'$ disjoint from $S$, and* *assume that the corresponding embedding problem is solvable. Then there exists a solution $L/K$ such that $L/k$ is only ramified in $S' \cup \{\mathfrak{q}\} \cup \{\mathfrak{p} \in \mathbb{P}(k) \mid \mathfrak{p}|\ell\}$, where $\mathfrak{q} \in \mathbb{P}(k) \setminus S$.*

*Proof.* By Kummer theory, any solution is of the form $N = K(\sqrt[\ell]{\alpha})$ for some $\alpha \in K^\times$. Let $(\alpha) = \prod_{i=1}^{r} \mathfrak{p}_i^{e_i}$ be the factorization of the principal ideal $(\alpha)$ in $K$, ordered such that $\mathfrak{p}_1, \ldots, \mathfrak{p}_s$ are unramified in $K/k$. Let $\tilde{e}_i \in \{0, \ldots, \ell-1\}$ such that $e_i \equiv \tilde{e}_i \pmod{\ell}$, $1 \leq i \leq s$. We claim that $\prod_{i=1}^{s} \mathfrak{p}_i^{\tilde{e}_i}$ is $\mathrm{Gal}(K/k)$-invariant, hence an ideal in $\mathcal{O}_k$. Indeed, let $\zeta \in k$ be a primitive $\ell$th root of unity. Then $\zeta^\sigma = \zeta$ for all $\sigma \in \mathrm{Gal}(K/k)$. Since the considered embedding problem is central, $N/K^{\langle \sigma \rangle}$ is abelian, so Lemma 2.8 implies that $\alpha^\sigma \equiv \alpha \pmod{K^{\times \ell}}$ for all $\sigma \in \mathrm{Gal}(K/k)$, as claimed.

Since in every ideal class there exists infinitely many prime ideals we can choose $\mathfrak{q} \in \mathbb{P}(k) \setminus S$ such that $\prod_{i=1}^{s} \mathfrak{p}_i^{\ell - \tilde{e}_i} \mathfrak{q}$ is a principal ideal of $\mathcal{O}_k$, generated by $b$ say. By [12], Thm. IV.7.2, the field $K(\sqrt[\ell]{b\alpha})$ is another solution to the embedding problem. By Theorem 3.6 it is at most ramified in $S' \cup \{\mathfrak{q}\} \cup \{\mathfrak{p} \in \mathbb{P}(k) \mid p|\ell\}$. □

We do not know whether in the situation described above there always exists a solution only ramified in $S' \cup \{\mathfrak{q}\}$. The following example shows that this is certainly not the case if the embedding problem is not central:

**Example 3.8.** Let $k = \mathbb{Q}(\sqrt{5})$ and consider the quadratic extension $K = k(\sqrt{\alpha})$ where $\alpha := 2 + \sqrt{5}$. Then for any $b \in \mathbb{Q}$ the extension $k(\sqrt{b\alpha})$ is ramified in 2.

**Corollary 3.9.** *In the situation of Proposition 3.3 assume that moreover $K/k$ is ramified in the finite set $S \subset \mathbb{P}(k)$. Then there exists a solution $L/k$ unramified outside $\{\mathfrak{q}\} \cup S \cup \{\mathfrak{p} \in \mathbb{P}(k) \mid \mathfrak{p}|\ell\}$ for some prime ideal $\mathfrak{q} \in \mathbb{P}(k)$.*

*Proof.* By Lemma 3.7 there exists an ideal $\tilde{\mathfrak{q}}$ and a solution $\tilde{L}/\tilde{k}$ to the induced embedding problem for $\tilde{K}/\tilde{k}$ with the desired properties. By the Theorem of Kochendörffer [12], Thm. IV.8.2, a solution field $L/k$ to the original embedding problem can be found inside the Galois closure of $\tilde{L}/k$. Clearly the latter is at most ramified in $\{\mathfrak{q}\} \cup S \cup \{\mathfrak{p} \in \mathbb{P}(k) \mid \mathfrak{p}|\ell\}$, where $\mathfrak{q}$ is the prime ideal of $\mathcal{O}_k$ lying under $\tilde{\mathfrak{q}}$. □

## 4. Lower bounds

We now obtain a lower bound for the rate of growth of the counting function for Galois extensions solving the embedding problem (*), using the corresponding result of Wright for cyclic groups of prime order, that is, the case $H = 1$.

**Theorem 4.1.** *Let $k$ be a number field, $H$ a finite group, $\ell$ a prime and*

$$1 \longrightarrow C_\ell \longrightarrow G \longrightarrow H \longrightarrow 1$$

*a central extension. Suppose that there exists a normal extension $L/k$ with Galois group $G$, linearly disjoint from $\tilde{k} = k(\zeta)$ over $k$. Then there exists $c(k, G) > 0$ such that*

$$c(k, G) \, x^{\frac{\ell}{(\ell-1)|G|}} \leq Z(k, G; x)$$

*for all $x \gg 0$.*

*Proof.* Denote by $K$ the fixed field in $L$ under the $C_\ell$-kernel. Then we are in the situation of Proposition 3.4. By Proposition 3.5(a) the correspondence between cyclic extensions of $k$ with group $C_\ell$ and the desired fields has finite fibers of bounded size. Moreover by 3.5(b) the discriminants are related by the given power relation up to bounded constants. Since

$$Z(k, C_\ell; x) \geq c_1(k, \ell) \, x^{a(C_\ell)} = c_1(k, \ell) \, x^{\frac{1}{\ell-1}}$$

for all large enough $x$ by the result of Wright [17], this implies that
$$Z(k, G; x) \geq c(k, \ell)\, x^{a(C_\ell)/|H|} = c(k, \ell)\, x^{\frac{\ell}{(\ell-1)|G|}}$$
for all $x \gg 0$. □

**Corollary 4.2.** *Let $k$ be a number field, $H$ a finite group, $\ell$ a prime such that $|H|$ is not divisible by primes smaller than $\ell$ and*
$$1 \longrightarrow C_\ell \longrightarrow G \longrightarrow H \longrightarrow 1$$
*a central extension. Let $S \subset \mathbb{P}(k)$ be finite and suppose that there exists a normal extension $L/k$ with Galois group $G$ not ramified in $S$. Then there exists $c(k, S, G) > 0$ such that*
$$c(k, S, G)\, x^{a(G)} \leq Z(k, S, G; x)$$
*for all $x \gg 0$.*

*Proof.* Note that by assumption $L$ is linearly disjoint over $k$ from the field $\tilde{k} = k(\zeta)$ obtained by adjoining the $\ell$th roots of unity. Moreover, by Remark 2.9 we have $a(G) = \frac{\ell}{(\ell-1)|G|}$ for the group $G$ in its regular permutation representation. So for $S = \emptyset$ the claim follows from Theorem 4.1.

In the general case, let $\alpha \in \tilde{K}^\times$ such that $\tilde{L} = \tilde{K}(\sqrt[\ell]{\alpha})$ as in Proposition 3.3. Let $b \in \tilde{k}^\times$ such that $k_b/k$ is a $C_\ell$-extension unramified in $S$. Then $\tilde{K}(\sqrt[\ell]{b\alpha})$ is a subfield of $\tilde{L}k_b$ and therefore unramified in $S \setminus \{\mathfrak{p} \in \mathbb{P}(k) \mid \mathfrak{p} | \ell\}$. Suppose $\mathfrak{p} \in S$ with $\mathfrak{p} \mid \ell$. Since $\tilde{L}k_b = L\tilde{k}_b$ and $\mathfrak{p}$ is unramified in $Lk_b/k$ we get that $\mathfrak{p}$ has ramification index dividing $\ell - 1$ in $\tilde{L}k_b/k$. By assumption $\ell - 1$ is coprime to $|G|$ and this implies that $\mathfrak{p}$ is unramified in $L_b/k$. Hence for any $b$ as above, the solution $L_b/k$ is unramified in $S$. The latter elements parametrize cyclic extensions of $k$ unramified in $S$. By Theorem 2.6 we know that their number can be estimated from below by
$$Z(k, S, C_\ell; x) \geq c_1\, x^{a(C_\ell)} = c_1\, x^{\frac{1}{\ell-1}}$$
for all sufficiently large $x$. We may now complete the argument as in the proof of Theorem 4.1. □

## 5. The upper bound

We turn to the problem of deriving an upper bound for the number of solutions to (*). Here we have to take into account all intermediate fields $K$.

**Theorem 5.1.** *Let $H$ be a finite group, $\ell$ a prime such that $|H|$ is not divisible by primes smaller than $\ell$ and*
$$1 \longrightarrow C_\ell \longrightarrow G \longrightarrow H \longrightarrow 1$$
*a central extension. Let $m \geq 1$ and assume that for any number field $k$ of degree $m = [k : \mathbb{Q}]$ and all $\delta > 0$ there exist constants $c(m, H, \delta)$, $c(k)$ such that $Z(k, H; x) \leq c(m, H, \delta)c(k)\, x^{a(H)+\delta}$. Then for all $\epsilon > 0$ there exists a constant $c(m, G, \epsilon)$ such that for all number fields $k$ of degree $m = [k : \mathbb{Q}]$ we have*
$$Z(k, G; x) \leq c(m, G, \epsilon)c(k)|\mathrm{Cl}(\tilde{k})_\ell|\, x^{a(G)+\epsilon} \qquad \text{for all } x > 0.$$

*Proof.* Let's fix a number field $k$ of degree $m = [k : \mathbb{Q}]$. Let $L/k$ be a Galois extension with Galois group isomorphic to $G$ and let $K$ be the fixed field of the central normal $C_\ell$-subgroup. According to Lemma 3.2 the field $L$ is now a solution to one of finitely many embedding problems, the number of which only depends on $H$ and $G$. Hence for any given $K$ it is sufficient to count solutions to a fixed embedding problem for the $H$-extension $K/k$.

Step 1: We compare solutions to the embedding problem above a fixed $K/k$ to cyclic extensions of $k$. So let $K/k$ be an $H$-extension and assume that the embedding problem for $K/k$ is solvable. Let $\mathfrak{q} \in \mathcal{O}_k$ be an ideal as in Corollary 3.9 and denote by $S_K \subset \mathbb{P}(k)$ the set of prime divisors of $\ell d_{K/k} \mathcal{O}_k$ union $\{\mathfrak{q}\}$. Then by Corollary 3.9 there exists a solution $L/k$ unramified outside $S_K$. Let $\alpha \in \tilde{K}$ such that $\tilde{L} = \tilde{K}(\sqrt[\ell]{\alpha})$.

Any other solution is of the form $L_b/k$, where $L_b$ is the solution field inside $\tilde{K}(\sqrt[\ell]{b\alpha})$, and where $b \in \tilde{k}^\times$ is such that $k_b/k$ is a $C_\ell$-extension of $k$. Moreover, by Proposition 3.5(a), once $\alpha$ is chosen the extension $k_b/k$ is uniquely determined by $L_b$. For $b \in \tilde{k}^\times$ let $d_b = \mathcal{N}(d_{k_b/k})$ be the norm of the discriminant of the $C_\ell$-extension $k_b/k$. Let $n := |H| = [K : k]$. The discriminant composition formula now shows that $d_{L_b/k} = d_{K/k}^\ell \mathcal{N}_{K/k}(d_{L_b/K})$. Moreover by the choice of $S_K$ we have $\mathcal{N}(d_{L_b/K})^{S_K} = (\mathcal{N}(d_{k_b/k})^{S_K})^n$, so

$$\mathcal{N}(d_{L_b/k}) \geq \mathcal{N}(d_{K/k})^\ell \left(\mathcal{N}(d_{k_b/k})^{S_K}\right)^n.$$

Denote by $Z_K(k, G; x)$ the number of solutions above $K$, that is, the number of $G$-extensions of $k$ such that $K$ is the fixed field of the $C_\ell$-kernel. Then the previous calculations show that

$$Z_K(k, G; x) \leq Y^{S_K}\left(k, C_\ell; \left(\frac{x}{D_K^\ell}\right)^{\frac{1}{n}}\right)$$

where we have set $D_K := \mathcal{N}(d_{K/k})$.

Step 2: We sum over all $H$-extensions $K/k$. We write

$$\mathfrak{K}(x) := \{K/k \text{ normal with group } H \text{ and } D_K^\ell \leq x\}$$

for the set of $H$-extensions of $k$ with bounded norm of the discriminant. Let $\epsilon > 0$ and $D_0 := \mathcal{N}(\mathfrak{a})(\leq c_1(m) d_k^{(\ell-1)/2} \leq c_1(m) d_K^{(\ell-1)/2})$. Then using our Theorem 2.5 to count cyclic extensions we get

$$\begin{aligned}
Z(k, G; x) &= \sum_{K \in \mathfrak{K}(x)} Z_K(k, G; x) \\
&\leq \sum_{K \in \mathfrak{K}(x)} Y^{S_K}\left(k, C_\ell; \left(\frac{x}{D_K^\ell}\right)^{\frac{1}{n}}\right) \\
&\leq \sum_{K \in \mathfrak{K}(x)} c(m, \delta_1) |\mathrm{Cl}(\tilde{k})_\ell| D_K^{r\delta_1} \left(\frac{x}{D_K^\ell}\right)^{\frac{1}{n(\ell-1)} + \frac{\delta_1}{n}} \\
&= c(m, \delta_1) |\mathrm{Cl}(\tilde{k})_\ell|\, x^{\frac{1}{n(\ell-1)} + \frac{\delta_1}{n}} \sum_{K \in \mathfrak{K}(x)} \frac{1}{D_K^{\frac{\ell}{n(\ell-1)} + \frac{\delta_1 \ell}{n} - r\delta_1}} \\
&\leq c(m, \delta_1) |\mathrm{Cl}(\tilde{k})_\ell|\, x^{\frac{1}{n(\ell-1)} + \epsilon} \sum_{d=1}^{\lfloor \sqrt[\ell]{x} \rfloor} \frac{l(d)}{d^{\frac{\ell}{n(\ell-1)} + \ell\epsilon - r\delta_1}}
\end{aligned}$$

for any $\delta_1 < n\epsilon$, where $l(d)$ is the number of $H$-extensions $K/k$ with norm of discriminant $D_K = d$. By our assumption on divisors of $|H|$ together with Remark 2.9 we have $a(H) \leq \ell/n(\ell-1)$. So

$$\sum_{d=1}^y l(d) = Z(k, H; y) \leq c(m, H, \delta_2) c(k) y^{a(H) + \delta_2} \leq c(m, H, \delta_2) c(k) y^{\frac{\ell}{n(\ell-1)} + \delta_2}$$

for all $\delta_2 > 0$ by our assumption on the number of $H$-extensions of $k$. By Lemma 2.4 this implies that the Dirichlet series

$$\sum_{d \geq 1} \frac{l(d)}{d^{\frac{\ell}{n(\ell-1)} + \ell\epsilon - r\delta_1}}$$

converges whenever $\ell\epsilon - r\delta_1 - \delta_2 > 0$. Moreover, since the upper bound for the partial sums, divided by $c(k)$, only depends on $m, H, \delta_2$, the infinite sum can again be bounded above in these terms. We thus finally obtain

$$Z(k, G; x) \leq c(m, H, \delta_1, \delta_2) \, c(k) |\text{Cl}(\tilde{k})_\ell| \, x^{\frac{1}{n(\ell-1)} + \epsilon}$$

whenever $r\delta_1 + \delta_2 < \ell\epsilon$. The desired result follows by choosing $\delta_1, \delta_2$ sufficiently small. □

Note that if we only assume that $Z(k, H; x) \leq c(k, H, \delta) \, x^{a(H)+\delta}$ for a particular number field $k$ then the proof of Theorem 5.1 still gives the conclusion $Z(k, G; x) \leq c(k, G, \epsilon) \, x^{a(G)+\epsilon}$ for this same field $k$.

## 6. Nilpotent groups

Combining the results from the previous sections we may conclude that Conjecture 1.1 holds for nilpotent Galois extensions of any number field. For the lower bound (with restricted ramification) we need the deep result, proved by Šafarevič, that there exists at least one extension (with controlled ramification) for each nilpotent group:

**Theorem 6.1.** (Šafarevič) *Let $k$ be a number field and $G$ a finite solvable group. Then for all finite subsets $S \subseteq \mathbb{P}(k)$ there exists a normal extension $K/k$ with Galois group $G$ which is unramified in $S$.*

For $S = \emptyset$ this is Theorem 9.5.1 in [14]. All the ingredients needed for the general statement are already contained in [14], but unfortunately, the authors do not prove it explicitly. In order to prove the above theorem we need to generalize Theorem 9.5.11 in [14]. We use the notation of [14]. Let $p$ be a prime number and $G$ be a finite group. Then $\mathcal{F}(d)$ is the free pro-$p$-$G$ operator group of rank $d$ (see [14], p. 482) which is a free pro-$p$ group of rank $|G|d$ with a certain $G$-action. For $v = (i, j)$ with $i \geq j \geq 1$ we denote by $\mathcal{F}(d)^{(v)}$ the filtration of $\mathcal{F}(d)$ defined in [14], p. 481, which is a refinement of the descending $p$-central series. We need to prove the following generalization of Theorem 9.5.11 in [14].

**Theorem 6.2.** *Let $K/k$ be a finite Galois extension of the number field $k$ and let $\varphi : G_k \twoheadrightarrow \text{Gal}(K/k) = G$. Then for every prime number $p$, all $n \in \mathbb{N}$ and all $v = (i, j)$, the split embedding problem*

$$1 \longrightarrow \mathcal{F}(n)/\mathcal{F}(n)^{(v)} \longrightarrow \mathcal{F}(n)/\mathcal{F}(n)^{(v)} \rtimes G \longrightarrow G \longrightarrow 1$$

*has a proper solution $N_n^v/k$. Furthermore, let $S \subset \mathbb{P}(k)$ be a finite set of primes which are unramified in $K/k$. Then we can choose the solution in such a way that the following conditions are satisfied:*

(1) *All primes $\mathfrak{p}$ which are ramified in $K/k$ or lying above $p$ or $\infty$ are completely decomposed in $N_n^v/K$.*
(2) *If $\mathfrak{p}$ is ramified in $N_n^v/K$, then $\mathfrak{p} \notin S$, $\mathfrak{p}$ splits completely in $K/k$ and $N_{n,\mathfrak{p}}^v/k_\mathfrak{p}$ is a cyclic totally ramified extension of local fields.*

*Proof.* For $S = \emptyset$ this is a special case of Theorem 9.5.11 in [14]. Let $\hat{S} := \{\mathfrak{p} \in S \mid \mathfrak{p} \text{ splits completely in } K/k\}$. Since $\hat{S}$ is finite there exists a quadratic extension $L/k$ such that all primes in $\hat{S}$ are inert and all primes which are ramified in $K/k$ are split. Then clearly $\text{Gal}(LK/k) = C_2 \times G$. We consider the corresponding embedding problem

$$1 \longrightarrow \mathcal{F}(n)/\mathcal{F}(n)^{(v)} \longrightarrow \mathcal{F}(n)/\mathcal{F}(n)^{(v)} \rtimes (G \times C_2) \longrightarrow G \times C_2 \longrightarrow 1$$

for $LK/k$, where the action of the $C_2$-part is trivial, i.e. $\mathcal{F}(n)/\mathcal{F}(n)^{(v)} \rtimes (G \times C_2) \cong C_2 \times (\mathcal{F}(n)/\mathcal{F}(n)^{(v)} \rtimes G)$. Using Theorem 9.5.11 in [14] we obtain a proper solution $\tilde{N}_n^\nu$ which is unramified for all primes in $S$ since these primes are not totally split in $LK/k$. A solution to our original embedding problem where all primes in $S$ are unramified is then given by the fixed field in $\tilde{N}_n^\nu$ of the direct factor $C_2$. Additionally conditions (1)+(2) are fulfilled. □

In the nilpotent case, which is the only one we are going to use, Theorem 6.1 now follows since a nilpotent group is the direct product of its Sylow $p$-subgroups, and any $p$-group is the quotient of $\mathcal{F}(n)$ for all sufficiently large $n$.

In the general solvable case, the proof of Theorem 6.1 works analogously as in [14], p. 507, where we use the stronger result of Theorem 6.2.

With this we obtain a proof of Conjecture 1.1 for nilpotent groups in their regular permutation representation.

**Theorem 6.3.** *Let $G$ be a finite nilpotent group in its regular representation and $k$ a number field. Then for all $\epsilon > 0$ there exists a constant $c_2(k, G, \epsilon)$, and for all finite subsets $S \subset \mathbb{P}(k)$ there exists a constant $c_1(k, S, G) > 0$ such that*

$$c_1(k, S, G)\, x^{a(G)} \leq Z(k, S, G; x) < c_2(k, G, \epsilon)\, x^{a(G)+\epsilon}$$

*for all $x \gg 0$.*

*Proof.* The proof is by induction on the order of $G$, the induction basis being given by the trivial group. Since $G$ is nilpotent, it is the direct product of its Sylow subgroups. Let $\ell$ denote the smallest prime divisor of $|G|$. Let $Z$ be a central cyclic subgroup of the Sylow $\ell$-subgroup of $G$, hence of $G$, of order $\ell$, and write $H := G/Z$. By the Theorem 6.1 of Šafarevič there exists a Galois extension $L/k$ with group $G$, unramified in $S$. Moreover, since $\ell$ is the smallest prime divisor of $|G|$, we have that $|G|$ is prime to $\ell - 1$, so $L/k$ is linearly disjoint from $k(\zeta_\ell)/k$. We are thus in the situation of Corollary 4.2, and the lower bound follows.

By induction, $Z(k, H; x)$ grows at most like $c(k, H, \delta)\, x^{a(H)+\delta}$, hence we are also in the situation of Theorem 5.1, and the upper bound is proved as well. □

Note that the upper bound follows without use of Šafarevič's theorem. Also note that induction actually allows to obtain a more precise result: the constant in the upper bound for nilpotent $G$ of order $|G| = \ell^s$ is of the form $c(m, G, \epsilon)|\mathrm{Cl}(\tilde{k})_\ell|^s$ for all $k$ of degree $m$ over $\mathbb{Q}$. In view of the result for abelian groups one might speculate that this could be improved to $c(m, G, \epsilon)|\mathrm{Cl}(\tilde{k})|_\ell$, not depending on $s$.

## 7. Changing the representation

The upper bound actually holds for arbitrary transitive permutation representations of $\ell$-groups. To prove this, we investigate the connection between the $a$-function for two different permutation representations of a finite group:

**Lemma 7.1.** *Let $G$ be a finite group with two faithful transitive permutation representations $\phi_1 : G \hookrightarrow \mathfrak{S}_n$, $\phi_2 : G \hookrightarrow \mathfrak{S}_m$. Define $a_i := a(\phi_i(G))$ and $\mathrm{ind}_i := \mathrm{ind} \circ \phi_i$, $i = 1, 2$. Assume that $a_2 \mathrm{ind}_2(\sigma) \geq a_1 \mathrm{ind}_1(\sigma)$ for all $\sigma \in G$.*

(a) *If $Z(k, \phi_1(G); x) \leq c_1 x^{a_1+\epsilon}$ for all $x > 0$ then there exists $c_2 > 0$ with*

$$Z(k, \phi_2(G); x) \leq c_2 x^{a_2 + \frac{a_2}{a_1}\epsilon} \qquad \text{for all } x > 0.$$

*Moreover, the quotient $c_2/c_1$ is bounded for fixed degree $m = [k : \mathbb{Q}]$.*

(b) *If $Z(k, \phi_2(G); x) \geq c_1 x^{a_2}$ for all $x \gg 0$, then there exists $c_2 > 0$ with*

$$Z(k, \phi_1(G); x) \geq c_2 x^{a_1} \qquad \text{for all } x \gg 0.$$

*Proof.* Denote the one-point stabilizers in the permutation representations $\phi_i$ by $H_i$, $i = 1, 2$. Let $N/k$ be a Galois extension with group $G$ and let $K_i := N^{H_i}$, $i = 1, 2$, be the fixed fields. Note that at most $(G : H_i)$ fields with group $\phi_i(G)$ correspond to the same Galois closure $N$, thus as far as asymptotics is concerned, this ambiguity can be neglected. A prime divisor $\mathfrak{p} \in \mathbb{P}(k)$ which is wildly ramified in $N/k$ necessarily divides $|G|$. Thus, by Lemma 2.1 we need not (and won't) consider wildly ramified primes. Now assume that $\mathfrak{p} \in \mathbb{P}(k)$ is tamely ramified in $N/k$ and let $\sigma$ be an inertia group generator at $\mathfrak{p}$. Then the precise power of $\mathfrak{p}$ dividing $d_{K_i/k}$ is $(d_{K_i/k})_\mathfrak{p} = \mathfrak{p}^{\mathrm{ind}_i(\sigma)}$, $i = 1, 2$ (see [10], p.100 and Prop. 6.3.1, for example). By assumption we thus have

$$(d_{K_2/k})_\mathfrak{p} = \mathfrak{p}^{\mathrm{ind}_2(\sigma)} \geq \mathfrak{p}^{\frac{a_1 \mathrm{ind}_1(\sigma)}{a_2}} = (d_{K_1/k})_\mathfrak{p}^{\frac{a_1}{a_2}}.$$

Since this is true for all but finitely many primes $\mathfrak{p}$ we conclude that

$$Z(k, \phi_2(G); x) \leq c\, Z(k, \phi_1(G); x^{\frac{a_2}{a_1}}) \leq c_2 x^{a_2 + \frac{a_2}{a_1}\epsilon}$$

for all $x > 0$ in part (a), with suitable constants $c, c_2$. Moreover, the proof shows that $c_2/c_1$ is bounded only depending on $m = [k : \mathbb{Q}]$.

Similarly, we obtain part (b). □

**Example 7.2.** Let $G$ be an $\ell$-group, $\phi_1 : G \hookrightarrow \mathfrak{S}_{|G|}$ the regular permutation representation and $\phi_2 : G \hookrightarrow \mathfrak{S}_n$ some faithful transitive permutation representation. Let $\sigma \in G$ have order $m$. We claim that

$$\mathrm{ind}_2(\sigma) \geq \frac{\ell(m-1)}{m(\ell-1)a_2(G)}.$$

Indeed, with $m = \ell^r$, write $b_1$ for the number of $\ell$-cycles of $\tau := \sigma^{\ell^{r-1}}$ in the chosen degree $n$ representation of $G$ and $b_2$ for the number of fixed points. Then $\mathrm{ind}_2(\tau) = n - (b_1 + b_2) = b_1\ell + b_2 - (b_1 + b_2) = b_1(\ell - 1)$. On the other hand, $\sigma$ has precisely $b_1/\ell^{r-1}$ cycles of length $\ell^r$, so

$$\mathrm{ind}_2(\sigma) \geq n - (b_1/\ell^{r-1} + b_2) = \frac{b_1(\ell^r - 1)}{\ell^{r-1}} = \frac{\mathrm{ind}_2(\tau)(\ell^r - 1)}{\ell^{r-1}(\ell - 1)}.$$

The result follows since $\mathrm{ind}_2(\tau) \geq a_2(G)^{-1}$ by definition.

Since

$$a_1(G) = \frac{\ell}{|G|(\ell-1)}, \quad \mathrm{ind}_1(\sigma) = \frac{|G|(m-1)}{m},$$

this shows that

$$\mathrm{ind}_2(\sigma) \geq \frac{a_1(G)\mathrm{ind}_1(\sigma)}{a_2(G)} \quad \text{for all } \sigma \in G.$$

**Corollary 7.3.** *Let $k$ be a number field, $\ell$ a prime, $G \leq \mathfrak{S}_n$ a transitive $\ell$-group, not necessarily in the regular representation. Then for all $\epsilon > 0$ there exists a constant $c(k, G, \epsilon) > 0$ such that*

$$Z(k, G; x) \leq c(k, G, \epsilon)\, x^{a(G)+\epsilon} \quad \text{for all } x > 0.$$

*Proof.* We reduce this to the case of the regular representation. By the previous example, the general assumption of Lemma 7.1 is satisfied. Moreover, by Theorem 6.3 applied to $G$ the upper bound is known to hold for the regular representation. Hence the result follows by Lemma 7.1(a). □

The above reasoning goes through whenever the smallest index in some permutation representation is attained on an element of smallest prime order. This is not always the case, as the next example shows.

**Example 7.4.** The smallest nilpotent group $G$ having a non-regular representation to which Lemma 7.1 is not applicable is the direct product of the cyclic group of order 2 with the non-abelian transitive subgroup of $\mathfrak{S}_9$ of order 27 and exponent 3. Indeed, for the regular representation we find $a_1 = 1/27$, while for the degree 18 representation of the direct product we have $a_2 = 1/8$. But for suitable elements $\sigma \in G$ of order 3 the indices are $\text{ind}_1(\sigma) = 36$, $\text{ind}_2(\sigma) = 8$, which violates the condition in Lemma 7.1.

In some special situations the proof of Theorem 5.1 works for permutation groups $G$ which are not given in regular representation. Assume that there exists a subgroup $U \leq G$ such that $\text{Core}(G, U) \cong C_\ell$ is contained in the center of $G$ and $G_1 \lhd U$ of index $\ell$, where $G_1$ denotes the point stabilizer. Let $H$ be the image of the action of $G$ on $G/U$. Then $G$ is a central extension

$$1 \longrightarrow C_\ell \longrightarrow G \longrightarrow H \longrightarrow 1.$$

We get a similar result as in Proposition 3.5.

**Proposition 7.5.** *Let $G$ and $H$ be as above and assume that $\ell$ is the smallest prime dividing $|G|$. Let $K/k$ be an extension with group $H$ (corresponding to the given permutation representation). Assume that the embedding problem is solvable. Then the correspondence $k_b \mapsto L_b$ between cyclic Galois extensions of $k$ with group $C_\ell$ and fields $L_b$ containing $K$ with $\text{Gal}(L_b/k) = G$ has finite fibers of size bounded only in terms of $H, \ell$, and $m$.*

*Moreover the discriminants satisfy*

$$c_1(K)\mathcal{N}(d_{L_b/k}) \leq \mathcal{N}(d_{k_b/k})^{|H|} \leq c_2(K)\mathcal{N}(d_{L_b/k})$$

*for constants $c_1(K), c_2(K) > 0$ only depending on $K$ (and $\ell$).*

*Proof.* Denote by $N$ the normal closure of $K/k$ and let $L/K$ be a field with $\text{Gal}(L/k) = G$. Adjoining the $\ell$th roots of unity we have that $\tilde{L} = \tilde{K}(\sqrt[\ell]{\alpha})$ for some $\alpha \in \tilde{K}$. The normal closure of $\tilde{L}/\tilde{k}$ equals $\tilde{N}(\sqrt[\ell]{\alpha}) =: M$. The correspondence $k_b \mapsto M_b$ is established in Proposition 3.5. The result follows from the fact that all solutions correspond to elements $\tilde{L}(\sqrt[\ell]{b\alpha})$, where $b$ is as in Proposition 3.5. □

**Theorem 7.6.** *Let $G$ and $H$ be as above and assume that $\ell$ is the smallest prime dividing $|G|$. Let $K/k$ be an extension with Galois group $H$. Assume that for all $\delta > 0$ there exists a constant $c(k, H, \delta)$ such that $Z(k, H; x) \leq c(k, H, \delta)\, x^{a(H)+\delta}$. Then for all $\epsilon > 0$ there exists a constant $c(k, G, \epsilon)$ such that*

$$Z(k, G; x) \leq c(k, G, \epsilon)\, x^{\frac{a(H)}{\ell}+\epsilon} \qquad \text{for all } x > 0.$$

*Proof.* Using Proposition 7.5 the proof is analogous to the proof of Theorem 5.1. □

**Example 7.7.** The group $G = \text{SL}_2(3)$ is a (non-split) central extension of $C_2$ by $H = \mathfrak{A}_4$. In its action on Sylow 3-subgroups, $G$ becomes a subgroup of $\mathfrak{S}_8$. The preceding Theorem shows that in this permutation representation of degree 8 we have $Z(k, \text{SL}_2(3); x) \leq c(k, G, \epsilon)\, x^{\frac{a(H)}{2}+\epsilon}$. Since $a(G) = \frac{1}{4} = \frac{a(H)}{2}$ this gives the upper bound predicted by Conjecture 1.1.

In the direction of lower bounds, we only obtain a partial result. For example it may happen that Lemma 7.1 gives the desired conclusion:

**Example 7.8.** (Groups of degree 8) In the following we denote by $T_i$ the $i$-th transitive group of degree $n$, where the degree is clear from the context. This is the group one gets by typing *TransitiveGroup(n,i);* in the computer algebra systems Gap or Magma.

The assumption of Lemma 7.1(b) is satisfied for the following transitive permutation groups of degree 8: The lower bound for $T_{14}$ (the symmetric group $\mathfrak{S}_4$ on eight letters) follows from the lower bound for this group on four letters proved in Proposition 10.2 (for the case $k = \mathbb{Q}$ it was shown to hold by Baily [1]). The lower bounds for $T_{28}$, respectively $T_{31}$, hold if they hold for $T_{27}$, resp. $T_{29}$.

But it seems that in order to derive the lower bound for arbitrary permutation representations we would need to show a much stronger result in the regular case, asserting some form of equi-distribution of inertia group generators. At least we may generalize Theorem 4.1 as follows:

**Proposition 7.9.** *Let $k$ be a number field and $1 \to C_\ell \to G \to H \to 1$ a central extension of finite groups. Let $\phi : G \hookrightarrow \mathfrak{S}_n$ be a faithful transitive permutation representation of $G$. Assume that there exists at least one Galois extension of $k$ with group $G$, linearly disjoint from $\tilde{k}/k$. Then there exists a constant $c(k,G) > 0$ such that*
$$c(k,G)\, x^{\frac{\ell}{(\ell-1)n}} \leq Z(k,G;x) \qquad \text{for all } x \gg 0.$$

*Proof.* Denote by $G_1$ a one-point stabilizer in $G$ in its given permutation representation $\phi$. Then $G_1$ cannot contain the central subgroup $C_\ell$, because otherwise $C_\ell$ would lie in the kernel of $\phi$. If $L/k$ is a Galois extension with group $G$, then the fixed field $N := L^{G_1}$ has group $\phi(G)$. The discriminants are related by the discriminant composition formula $d_{L/k} = d_{N/k}^{|G_1|}\mathcal{N}(d_{L/N})$. The second factor $\mathcal{N}(d_{L/N})$ is bounded in terms of the discriminant of the fixed field $K := L^{C_\ell}$. Thus if we consider extensions $L$ above a chosen field $K/k$ with group $H$, then $d_{L/k}$ and $d_{N/k}^{|G_1|}$ only differ by bounded terms. Now the proof of Theorem 4.1 shows that the number of $G$-extensions above a fixed $K$ grows at least like $c\, x^{\frac{\ell}{(\ell-1)|G|}}$, so our assertion follows from $|G : G_1| = n$.  $\square$

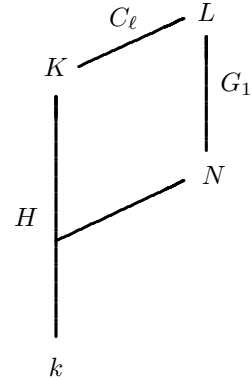

For an inductive approach, the lower bound would have to depend on $G$ in some more precise way. Still, in some cases this gives the correct lower bound even for non-regular permutation representations:

**Example 7.10.** We continue Example 7.7 for the group $G = \mathrm{SL}_2(3)$ in its degree 8 representation. It is well-known that there exist Galois extensions of any number field $k$ with group $G$ (see e.g. [12], Appendix, Table 6). Applying Proposition 7.9 we conclude that
$$Z(k, \mathrm{SL}_2(3); x) \geq c(k)\, x^{\frac{1}{4}}$$
for some $c(k) > 0$. It is easily seen that $a(G) = 1/4$ in this representation, so together with Example 7.7 this proves Conjecture 1.1 for $\mathrm{SL}_2(3) \leq \mathfrak{S}_8$.

## 8. Wreath products

At the moment, we can prove the lower bound for wreath products only when the kernel is cyclic of prime order. Let $\ell$ be a prime number.

**Lemma 8.1.** *Let $K$ be a number field and $S \subset \mathbb{P}(K)$ finite. Then there exists a prime $p \neq \ell$ and a cyclic extension $M/K$ unramified in $S$ such that*
  (1) *$p$ is totally split in $K$,*
  (2) *there exists exactly one prime ideal $\mathfrak{p} \subseteq \mathcal{O}_K$ above $p$ such that $\mathfrak{p}$ is ramified in $M$,*
  (3) *no divisor of $p$ lies in $S$.*

*Proof.* By [13], Cor. 7 to Prop. 7.9, there exist infinitely many primes $p \equiv 1 \pmod{\ell}$ which split totally in $K$. Denote the unit rank of $K$ by $r$ and let $p_1, \ldots, p_{r+2}$ be primes with the above property, such that moreover no divisor of $p_i$ lies in $S$. For $1 \leq i \leq r+2$, choose one prime ideal $\mathfrak{p}_i$ of $\mathcal{O}_K$ lying above $p_i$. Let $\mathfrak{a} := \mathfrak{p}_1 \cdots \mathfrak{p}_{r+2}$.

We get that $C_\ell^{r+2}$ is a factor group of $(\mathcal{O}/\mathfrak{a})^\times$. Denote by $A$ the ray class group of $\mathfrak{a}$. Then by [9], pp. 126–127, we have

$$A/\mathrm{Cl}(K) \cong (\mathcal{O}/\mathfrak{a})^\times/(U_K \cap K_\mathfrak{a}),$$

where $U_K$ denotes the unit group of $\mathcal{O}_K$ and $K_\mathfrak{a}$ denotes the subgroup of $K^\times$ of elements congruent 1 (mod $\mathfrak{a}$). Since $U_K/U_K^\ell$ has at most $\ell$-rank $r+1$, $C_\ell$ must be a factor group of $A/\mathrm{Cl}(K)$. Therefore there exists a $C_\ell$-extension $M/K$ which is at most ramified in the ideals dividing $\mathfrak{a}$. Furthermore, by construction it must be ramified in at least one divisor $\mathfrak{p}_i$ of $\mathfrak{a}$. $\square$

**Theorem 8.2.** *Let $K/k$ be an extension of number fields with Galois group $H$ and $\ell$ a prime number. Then for any finite $S \subset \mathbb{P}(k)$ above which $K/k$ is unramified there exists a constant $c = c(k, S, \ell, H) > 0$ such that*

$$c\,x^{\frac{1}{\ell-1}} \leq Z(k, S, C_\ell \wr H; x)$$

*for all $x \gg 0$.*

*Proof.* Let $M$ be a $C_\ell$-extension of $K$ as given in Lemma 8.1, not ramified above $S$, and with corresponding prime $p \neq \ell$. Let $S'$ be the set of prime divisors of $d_{M/K}$ in $\mathcal{O}_K$ together with the primes above $S$. By Theorem 2.6 we have $Z(K, S', C_\ell; x) \geq d\,x^{\frac{1}{\ell-1}}$, for sufficiently large $x$, for some constant $d > 0$ depending on $K, S', \ell$ only.

Let $N$ be such a $C_\ell$-extension of $K$ unramified in $S'$. Then $MN/K$ contains a subfield $\tilde{N}$ of degree $\ell$ different from both $M$ and $N$. Since $M, N$ are ramified in relatively coprime ideals, the discriminant of $\tilde{N}/K$ is given by $d_{\tilde{N}/K} = d_{M/K} d_{N/K}$. Furthermore, by construction $\tilde{N}/k$ is unramified above $S$, and exactly one prime ideal of $K$ above $p$ ramifies in $\tilde{N}$. In particular, $\tilde{N}/K$ is distinct from all its $K/k$-conjugates, and $\tilde{N}/k$ has Galois group $C_\ell \wr H$.

Thus any choice of $N$ leads to a $C_\ell \wr H$-extension of $k$ unramified in $S$. Suppose $\tilde{N}_1$ and $\tilde{N}_2$ are isomorphic. Then the corresponding fields $N_1$ and $N_2$ must be contained in $M\tilde{N}_1$. If $N_1 \neq N_2$ then $M\tilde{N}_1 = N_1 N_2$. The latter extension is unramified in $p$ which is a contradiction. Therefore their number can be bounded below by

$$Z(k, S, C_\ell \wr H; x) \geq \frac{d}{r}\left(\frac{x}{|d_M|}\right)^{\frac{1}{\ell-1}} = \frac{d}{r|d_M|^{\frac{1}{\ell-1}}} x^{\frac{1}{\ell-1}}$$

for all sufficiently large $x$, where $r$ is the number of non-isomorphic extensions $\tilde{N}/K$ which are isomorphic over $k$. Clearly the constant only depends on $H$ and $\ell$, and the proof is complete. $\square$

Using a general result from [11] we obtain the correct upper bound for wreath products of $C_2$ with arbitrary groups:

**Proposition 8.3.** *Let $k$ be a number field, $H \leq \mathfrak{S}_n$ a transitive permutation group such that $Z(k, H; x) \leq c(k, H, \delta)\, x^{a(H)+\delta}$ for all $\delta > 0$. Then for any $\epsilon > 0$ there exists a constant $c(k, 2 \wr H, \epsilon)$ such that*

$$Z(k, 2 \wr H; x) \leq c(k, 2 \wr H, \epsilon)\, x^{a(2 \wr H)+\epsilon}.$$

*Proof.* The result is certainly true if $k$ has no $H$-extensions. So now assume that $K/k$ is a $H$-extension. By Theorem 2.5 the number of $C_2$-extensions of any number field $K$ of degree $m = [K:\mathbb{Q}]$ grows at most like $c(m,\epsilon)\,|\mathrm{Cl}(K)_2|\,x^{1+\epsilon}$. But the class number and hence also the 2-part of the class group grows at most like $c(\delta)\,|d_K|^{0.5+\delta}$ for any $\delta > 0$ (see [13], Thm 4.4). Thus we are in the situation of Cor. 5.3 in [11], and the result follows. $\square$

Theorem 8.2 together with Proposition 8.3 show that Conjecture 1.1 is consistent with taking wreath products with the cyclic group $C_2$:

**Corollary 8.4.** *Let $k$ be a number field and $H \leq \mathfrak{S}_n$ a transitive permutation group for which*
$$1 \leq Z(k, H; x) \leq c(k, H, \epsilon) \, x^{a(H)+\epsilon}$$
*for all $\epsilon > 0$ and all $x \gg 0$. Then with $G := C_2 \wr H$ there exist for all $\epsilon > 0$ constants $c_1(k, G), c_2(k, G, \epsilon) > 0$ such that we have*
$$c_1(k, G) \, x^{a(G)} \leq Z(k, G; x) \leq c_2(k, G, \epsilon) \, x^{a(G)+\epsilon}$$
*for all $x \gg 0$.*

## 9. Direct products

Let $H, Z$ be transitive permutation groups of degrees $n$, $m$ respectively. Then the direct product $G := H \times Z$ has a natural transitive permutation representation on $nm$ points. It was shown in [11], Lemma 4.1, that with respect to these representations we have
$$a(G) = \max\{\frac{a(H)}{m}, \frac{a(Z)}{n}\}.$$

We apply the results of the previous sections to investigate direct products with nilpotent groups.

**Theorem 9.1.** *Let $H$ be a transitive permutation group of degree $n$, $Z$ a regular nilpotent permutation group of order $m$ and $G := H \times Z$. Assume that either*

(1) $a(H) \leq \frac{ma(Z)}{n}$, or
(2) *for any finite $S \subset \mathbb{P}(k)$ there exists a constant $c(S, H) > 0$ such that we have $c(S, H) \, x^{a(H)} \leq Z(k, S, H; x)$ for all large enough $x$.*

*Then for any finite $S \subset \mathbb{P}(k)$ for which there exists an $H$-extension of $k$ unramified in $S$, there exists a constant $c(S) > 0$ such that*
$$c(S) \, x^{a(G)} \leq Z(k, S, G; x)$$
*for all $x \gg 0$.*

*Proof.* By Corollary 4.2 the number of $Z$-extensions $L$ of $k$ unramified in any given finite subset $S \subset \mathbb{P}(k)$ grows at least like $c'(S) \, x^{a(Z)}$ for some $c'(S) > 0$. In the first case we have $a(H) \leq ma(Z)/n$, so the assertion follows from the general result for direct products in [11], Prop. 4.2. Similarly, in the second case we may assume that $a(H) > ma(Z)/n$, and by the assumption on $H$-extensions the assertion is covered by the same general result. □

## 10. Examples

In this section we investigate which transitive permutation groups of low degree can be handled by the methods developed in this paper.

In degree 3 we have the following result of Datskovsky and Wright [5] for the symmetric group on three letters. We denote by $\zeta_k(1)$ the residue at 1 of the Dedekind $\zeta$-function.

**Theorem 10.1.** (Datskovsky and Wright) *Let $k$ be a number field with signature $(r_1, r_2)$. Then*
$$Z(k, S_3; x) \sim \frac{2}{3}^{r_1-1} \frac{1}{6}^{r_2} \frac{\zeta_k(1)}{\zeta_k(3)} x.$$

In degree 4, there exist three non-abelian transitive groups. Extensions with the dihedral group $D_4$ were counted by Cohen, Diaz y Diaz and Olivier [3] whose (more precise) result confirms Conjecture 1.1. Baily [1], Thm. 1 and 3, has proved the following lemma in the case $k = \mathbb{Q}$.

**Proposition 10.2.** *The lower bound bound predicted by Conjecture 1.1 for the alternating and symmetric group on four letters is correct.*

*Proof.* Let $M/k$ be a degree 3 extension of number fields with Galois group $C_3$ or $\mathfrak{S}_3$, respectively. Let $L = M(\sqrt{\alpha})$, where $\alpha \in M$ is squarefree. In [7], Lemma 12, it is shown that $L/k$ has Galois group $6T4 = A_4(6)$ or $6T7 = S_4(6)$, respectively, if and only if $\mathcal{N}_{M/k}(\alpha)$ is a square. Let $\epsilon$ be a fundamental unit of $M$ with norm $a := \mathcal{N}_{M/k}(\epsilon)$. Then $a\epsilon$ is squarefree in $M$ and $\mathcal{N}_{M/k}(a\epsilon) = a^3 \mathcal{N}_{M/k}(\epsilon) = a^4$. The extension $L = M(\sqrt{a\epsilon})$ is at most ramified in prime ideals above 2. Denote the corresponding degree 4 extension in the normal closure of $L/k$ by $K$. Then in [7], Cor. 6, it is proved that $d_K = d_M \mathcal{N}(d_{L/M})$ holds. As in the proof of Lemma 2.1 we can bound $\mathcal{N}(d_{L/M})$ using a constant depending on $m = [k : \mathbb{Q}]$. Thus for each $\mathfrak{A}_3$-extension respectively $\mathfrak{S}_3$-extension of $k$ we have constructed an $\mathfrak{A}_4$-extension respectively $\mathfrak{S}_4$-extension of $k$ with the the same 2'-part of the discriminant. These extensions are all different since the normal closures have different degree 3 subfields. Since $a(\mathfrak{A}_4) = \frac{1}{2} = a(\mathfrak{A}_3)$ and $a(\mathfrak{S}_4) = 1 = a(\mathfrak{S}_3)$ the results on lower bounds for the cyclic group $\mathfrak{A}_3$ respectively for $\mathfrak{S}_3$ in Theorem 10.1 allow to conclude. □

In degree 6, there exist 16 transitive permutation groups up to equivalence, twelve of which are solvable. Our results allow to determine the growth of $Z(k, G; x)$ for the following groups:

**Proposition 10.3.** *Conjecture 1.1 holds for the transitive groups*

$$T_i \quad \text{with } i \in \{1, 2, 3, 6, 11\},$$

*of degree 6.*

*Proof.* The groups $T_1, T_2$ are abelian, the group $T_3$ is the direct product $\mathfrak{S}_3 \times C_2$, so Theorem. 7.6 together with Theorem 9.1 applies. $T_6$ and $T_{11}$ are wreath products $2 \wr 3$ and $2 \wr \mathfrak{S}_3$, hence we may apply Cor. 8.4. □

Moreover, partial results which follow from our theorems are collected in the subsequent table:

Bounds for degree 6 extensions of $k$

| Nr. | $|G|$ | $G$ | lower | $a(G)$ | upper | remarks |
|---|---|---|---|---|---|---|
| 4 | 12 | $\mathfrak{A}_4$ |  | 1/2 | 1/2 | Lemma 7.1 + [1] |
| 5 | 18 | $3 \wr 2$ | 1/2 | 1/2 |  | Thm. 8.2 |
| 7 | 24 | $\mathfrak{S}_4$ | 1/2 | 1/2 |  | Lemma 7.1 + Prop. 10.2 |

In degree 8, there exist 50 transitive permutation groups up to equivalence, 45 of which are solvable. Our results allow to determine the growth of $Z(k, G; x)$ for the following groups:

**Proposition 10.4.** *Conjecture 1.1 holds for the transitive groups*

$$T_i \quad \text{with } i \in \{1, 2, 3, 4, 5, 9, 12, 27, 28, 31, 35\},$$

*of degree 8.*

*Proof.* The groups $T_i$, $i = 1, 2, 3$ are abelian. The groups $T_4, T_5$ are regular nilpotent, so the claim follows from Theorem 6.3. The group $T_9$ is the direct product $C_2 \times D_4$, for which we may apply Cor. 7.3 together with Theorem 9.1 and [3]. The group $T_{12} \cong \mathrm{SL}_2(3)$ is treated in Examples 7.7 and 7.10. The groups $T_{27}, T_{31}, T_{35}$ are wreath products of $C_2$ with $C_4, C_2^2, D_4$ respectively, to which we apply Cor. 8.4. Finally, the lower bound for $T_{28}$ follows by application of Lemma 7.1 to the lower bound for $T_{27}$ (see Example 7.8), while the upper bound is Cor. 7.3. □

For the remaining groups, our results allow to determine lower respectively upper bounds for the growth of $Z(k, G; x)$ in the following cases:

Bounds for degree 8 extensions of $k$

| Nr. | $|G|$ | $G$ | lower | $a(G)$ | upper | remarks |
|---|---|---|---|---|---|---|
| 6 | 16 | $D_8$ | 1/4 | 1/3 | 1/3 | Prop. 7.9+Cor. 7.3 |
| 7 | 16 | | 1/4 | 1/2 | 1/2 | Prop. 7.9+Cor. 7.3 |
| 8 | 16 | | 1/4 | 1/3 | 1/3 | Prop. 7.9+Cor. 7.3 |
| 10 | 16 | | 1/4 | 1/2 | 1/2 | Prop. 7.9+Cor. 7.3 |
| 11 | 16 | | 1/4 | 1/2 | 1/2 | Prop. 7.9+Cor. 7.3 |
| 13 | 24 | $\mathfrak{A}_4 \times 2$ | 1/4 | 1/4 | | Thm. 9.1 |
| 14 | 24 | $\mathfrak{S}_4$ | 1/4 | 1/4 | | Lemma 7.1+Prop. 10.2 |
| 15 | 32 | | 1/4 | 1/2 | 1/2 | Prop. 7.9+Cor. 7.3 |
| 16 | 32 | | 1/4 | 1/2 | 1/2 | Prop. 7.9+Cor. 7.3 |
| 17 | 32 | $4 \wr 2$ | 1/4 | 1/2 | 1/2 | Prop. 7.9+Cor. 7.3 |
| 18 | 32 | $2^2 \wr 2$ | 1/4 | 1/2 | 1/2 | Prop. 7.9+Cor. 7.3 |
| 19 | 32 | | 1/4 | 1/2 | 1/2 | Prop. 7.9+Cor. 7.3 |
| 20 | 32 | | 1/4 | 1/2 | 1/2 | Prop. 7.9+Cor. 7.3 |
| 21 | 32 | | 1/4 | 1/2 | 1/2 | Prop. 7.9+Cor. 7.3 |
| 22 | 32 | | 1/4 | 1/2 | 1/2 | Prop. 7.9+Cor. 7.3 |
| 24 | 48 | $\mathfrak{S}_4 \times 2$ | 1/4 | 1/2 | | Prop. 7.9+Prop. 10.2 |
| 26 | 64 | | 1/4 | 1/2 | 1/2 | Prop. 7.9+Cor. 7.3 |
| 29 | 64 | | 1/4 | 1/2 | 1/2 | Prop. 7.9+Cor. 7.3 |
| 30 | 64 | | 1/4 | 1/2 | 1/2 | Prop. 7.9+Cor. 7.3 |
| 38 | 192 | $2 \wr \mathfrak{A}_4$ | 1 | 1 | | Cor. 8.4 |
| 44 | 384 | $2 \wr \mathfrak{S}_4$ | 1 | 1 | | Cor. 8.4 |

J.K.: Universität Heidelberg, IWR, Im Neuenheimer Feld 368, 69120 Heidelberg, Germany.
*E-mail address*: klueners@iwr.uni-heidelberg.de

G.M.: FB Mathematik/Informatik, Universität Kassel, Heinrich-Plett-Strasse 40, 34132 Kassel, Germany.
*E-mail address*: malle@mathematik.uni-kassel.de